\documentclass{amsart}
\usepackage{amssymb,latexsym}

\theoremstyle{plain}
\newtheorem{theorem}{Theorem}

\newtheorem{lemma}{Lemma}

\newtheorem{remark}{Remark}

\theoremstyle{definition}
\newtheorem{definition}{Definition}

\theoremstyle{remark}

\numberwithin{equation}{section}

\begin{document}
\title[completion and Caratheodory Extension ]{An isomorphism between the completion of an algebra and its Caratheodory Extension}
\author{Jun Tanaka}
\address{University of California, Riverside, USA}
\email{juntanaka@math.ucr.edu, yonigeninnin@gmail.com, junextension@hotmail.com}

\date{January, 10, 2008}

\begin{abstract}
  Let $\Omega$ denote an algebra of sets and $\mu$ a $\sigma$-finite measure. We then prove that the completion of $\Omega$ under the pseudometric $d(A,B)$ = $\mu^{\ast}(A \triangle B)$ is $\sigma$-algebra isomorphic and isometric to the Caratheodory Extension of $\Omega$ under the equivalence relation $\sim$.

\end{abstract}

\maketitle


\section{Introduction}
This paper shows a new result by combining two papers authored by P.F. Mcloughlin and myself ([4], [5]).

Let $\mu$ be a $\sigma$-finite measure and let $\Omega$ denote an algebra of sets; i.e., $\Omega$ is closed under unions and complements. Let ($X$, $\Omega$, $\mu$) denote a measure space where $\mu (X)$ is finite from the $\sigma$-finite property. Let $\mu^{\ast}$ denote the outer measure defined by $\mu^{\ast}(A)$  = $\inf \{ \sum \mu (A_{i}) \mid  E  \subseteq \cup A_{i}$ and $A_{i}  \in \Omega$ for all i $\geq $ 1$\}$, for any $A \in \mathbf{P}(X)$ where $\mathbf{P}(X)$ is the power set of $X$. Clearly, $d(A,B)$ = $\mu^{\ast}(A \triangle B)$ is a pseudometric, where $\triangle$ is the symmetric difference of sets. In addition, $d$ is a metric on $\mathbf{P}(X)_{\diagup_{\sim} } $, where A $\sim$ B iff $\mu^{\ast}(A \triangle B)$ =0. From [1], pg 292, $\mu^{\ast} |_{\Omega} = \mu $.

In [5], we defined a $\mu  $-Cauchy sequence $\{B_{n} \}$, $B_{n} \in \Omega$, if $\lim \mu(B_{n} \triangle B_{m})$ $\rightarrow$ 0 as $n,m \rightarrow \infty$. Let $ \widetilde{\mathbf{S}}$ = $\{S \in \mathbf{P}(X) \ |$ $ \exists \ \mu$-Cauchy sequence  $\{B_{n}   \}   $ s.t. $ \lim \mu^{\ast} ( B_{n}   \triangle S)  =0   \}$.

In the first joint paper [5], we proved that $ \widetilde{\mathbf{S}}$ is a $\sigma$-algebra where, for any $\mu$-Cauchy sequence $ \{ B_{n} \}$ such that $ \lim \mu^{\ast} ( B_{n}   \triangle S)  =0 $, the measure $\widetilde{\mu} (S) $ on $ \widetilde{\mathbf{S}}$ is defined as $\widetilde{\mu} (S) $ = $\lim \mu ( B_{n}   )   $. In addition, we proved that $  \widetilde{\mu}$ is a countably additive measure on $ \widetilde{\mathbf{S}}$. Thus, ($\widetilde{\mu}$, $ \widetilde{\mathbf{S}}$) is a measure space. We showed that the Caratheodory Extension of $\Omega$ can be expressed as the set of limit points of $\mu$-Cauchy sequences under the pseudometric $d(A,B)$ = $\mu^{\ast}(A \triangle B)$. Moreover, when the measure is a sigma finite measure, we obtained an equivalent expression of the Caratheodory Extension, $\{S \in \mathbf{P}(X) $ $ | $  $\exists \ \mu$-Cauchy sequence  $\{B_{n}   \}   $ s.t. $ \lim \mu^{\ast} ( B_{n}   \triangle S)  =0   \}$. Theorem 2 in [5] shows that $E$ is a measurable set iff  $E$ is in $ \widetilde{\mathbf{S}}$. Thus, the measure space ($\widetilde{\mu}$, $ \widetilde{\mathbf{S}}$) agrees with the Caratheodory Extension when $\mu$ is a finite measure. Moreover, it shows that measurable sets are exactly limit points of $\mu$-Cauchy sequences. The $\sigma$-finite case follows from the finite case.

 From the second joint paper [4], we denoted by ($\overline{d}$, $\overline{\Omega}$) the completion of ($d$, $\Omega_{\diagup_{\sim} } $). Let $\mathcal{S}$ be the set of all $\mu  $-Cauchy sequences in ($d$, $\Omega_{\diagup_{\sim} } $). By the completion procedures, we know $\{B^{\alpha}_{n} \}$ $\sim$ $\{B^{\gamma}_{n} \}$ iff $\lim d(B^{\alpha}_{n} ,B^{\gamma}_{n} )$ = 0 defines an equivalence relation on $\mathcal{S}$. Moreover, $\overline{\Omega}$ = $\mathcal{S}_{\diagup_{\sim} }$  and  $\overline{d}( \overline{\{B^{\alpha}_{n}\}}, \overline{\{B^{\gamma}_{n}\}} ) $ = $\lim d( \{B^{\alpha}_{n}\},\{B^{\gamma}_{n}\} ) $, where $\overline{\{B^{\alpha}_{n}\}}$ is the class of $\{B^{\alpha}_{n}\}$. Let $E_{\alpha}$ = $\overline{\{B^{\alpha}_{n}\}}$ and $E_{A}$ = $\overline{\{ A \}}$ when $A \in \Omega$. Let $\overline{\mu}(E_{\alpha})$ = $\overline{d}(E_{\alpha},E_{\emptyset})$ = $\lim d(B^{\alpha}_{n},\emptyset)$ = $\lim \mu(B^{\alpha}_{n})$. Note that in [4] $d(A,B)$ := $\mu(A \triangle B)$, whereas in [5] $d(A,B)$ := $\mu^{\ast}(A \triangle B)$; the completion of $\Omega$ will remain the same due to the property $\mu^{\ast} |_{\Omega} = \mu $. Note that $\overline{\mu}(E_{A})$ = $\mu(A)$ when $A \in \Omega$. In [4], we defined set-theoretic notations for unions, intersections, and complements on $\overline{\Omega}$ as follows: $\bigcup$ : $\overline{\Omega} \times  \overline{\Omega}  \rightarrow    \overline{\Omega} $ where $\bigcup(E_{\alpha}  \times  E_{\gamma} )$  =  $E_{\alpha} \bigcup  E_{\gamma}$ = $\overline{ \{B^{\alpha}_{n} \cup  B^{\gamma}_{n} \}}$; similarly for intersections on $\overline{\Omega}$. $\cdot^{\textbf{C}}  $: $\overline{\Omega}  \rightarrow    \overline{\Omega} $ where $(\overline{\{ B^{\alpha}_{n} \} }  )^{\textbf{C}} $ = $\overline{\{ (B^{\alpha}_{n})^{\textbf{C}} \} }$ and we showed the set theoretic notations are well defined on $\overline{\Omega}$ in [4].

 We showed the set theoretic notations are well defined on $\overline{\Omega}$ in [4]. Note that $E_{\alpha} \bigcap E_{\gamma}$ = $E_{\emptyset}$ iff $\overline{\mu}(E_{\alpha} \bigcap E_{\gamma})$ = 0 iff $\lim \mu (    B^{\alpha}_{n}   \cap  B^{\gamma}_{n}     )$ = 0.
We say  $E_{\alpha}$ and $E_{\gamma}$ are disjoint iff  $E_{\alpha} \bigcap E_{\gamma}$ = $E_{\emptyset}$. Thus, if $E_{\alpha_{1}}$ and $E_{\alpha_{2}}$ are disjoint, then $\overline{\mu}(E_{\alpha_{1}} \bigcup E_{\alpha_{2}})$ = $\overline{\mu}(E_{\alpha_{1}})$  +   $\overline{\mu}(E_{\alpha_{2}})$ .

As for the infinite union on $\overline{\Omega}$; if $E_{\alpha_{i}}$ $\in \overline{\Omega}$ for i $\geq$ 1, there exists a unique E := $\bigcup_{i=1}^{\infty}   E_{\alpha_{i}}  $ in $\overline{\Omega}$ such that $\bigcup_{i=1}^{n}   E_{\alpha_{i}}  \subset$ E for all n, and $\lim  \overline{\mu}(E \bigcap ( \bigcup_{i=1}^{n}     E_{\alpha_{i}}) ^{\textbf{C}}  )   $ = 0.

 In addition, we showed that for any $\mu$-Cauchy sequence $\{ B_{n} \}$, there exists a $f(n)$ $>$ n such that $ \lim \mu^{\ast} ( B_{n}   \triangle \overline{\lim} B_{f(n)}  )  =0   $.

In this paper, I define a $\sigma$-algebra isomorphism between two $\sigma$-algebras, and define a map $F :$ $\overline{\Omega}   \rightarrow \mathbf{P}(X) $ given by $F ( \overline{\{B_{n}\}} )  = \overline{\lim} B_{f(n)}  $ where $f(n)$ is defined as above. We will show that $F$ is an isometry and a $\sigma$-algebra isomorphism between the completion $\overline{\Omega}$ and the Caratheodory Extension of $\Omega$ under the equivalence relation $\sim$ defined as A $\sim$ B iff $\mu^{\ast}(A \triangle B)$ =0.






\section{Main Result}

\begin{definition}
For A, B in $\mathbf{P}(X) $, A = B a.e. iff $\mu^{\ast} ( A  \triangle B )  =0. $

\end{definition}

\begin{definition}\label{d:5}
Define a map $F :$ $\overline{\Omega}   \rightarrow \mathbf{P}(X) $ given by $F ( \overline{\{B_{n}\}} )  = \overline{\lim} B_{f(n)}  $ where $ \lim \mu^{\ast} ( B_{n}   \triangle \overline{\lim} B_{f(n)}  )  =0   $. Note that such $f(n)$ always exists by Lemma 20 in [4].

\end{definition}

\begin{remark}
$F$ is a map into $ \widetilde{\mathbf{S}}$ by the definition of $ \widetilde{\mathbf{S}}$.
\end{remark}

\begin{lemma}

$F$ is well defined.

\begin{proof}
Suppose that $\overline{\{A_{n}\}} =  \overline{\{B_{n}\}}$.

There exist f(n) and g(n) such that $ \lim \mu^{\ast} ( A_{n}   \triangle \overline{\lim} A_{f(n)}  )  =0   $  and  $ \lim \mu^{\ast} ( B_{n}   \triangle \overline{\lim} B_{g(n)}  )  =0   $ by Lemma 20 in [4].

$ \mu^{\ast} (  \overline{\lim} A_{f(n)}  \triangle \overline{\lim} B_{g(n)}  )  \leq  $ $ \mu^{\ast} (  \overline{\lim} A_{f(n)}  \triangle  A_{f(n)}  )  +   \mu^{\ast} (   A_{f(n)}  \triangle B_{g(n)}  ) +   \mu^{\ast} (  B_{g(n)}  \triangle \overline{\lim} B_{g(n)}  )  $ by the triangle inequality.

By taking the limit on both sides, $ \mu^{\ast} (  \overline{\lim} A_{f(n)}  \triangle \overline{\lim} B_{g(n)}  ) = 0$.

Thus, $  \overline{\lim} A_{f(n)} =  \overline{\lim} B_{g(n)}  $ a.e.. Therefore, F is well-defined.

\end{proof}
\end{lemma}

\begin{theorem}\label{t:3}
F is an isometry between $\overline{\Omega}$ and $ \widetilde{\mathbf{S}}_{\diagup_{\sim} }$.

\begin{proof}
First, we show F is onto $ \widetilde{\mathbf{S}}$. Let $X \in \widetilde{\mathbf{S}}$. Then there exists a $\mu$-Cauchy sequence  $\{B_{n}   \}   $ such that $ \lim \mu^{\ast} ( B_{n}   \triangle X)  =0   $.

Then there exist f(n) such that $ \lim \mu^{\ast} ( B_{n}   \triangle  \overline{\lim} B_{f(n)}   )  =0   $. Thus $F ( \overline{\{B_{n}\}} )  = \overline{\lim} B_{f(n)}  = X $ a.e.. Therefore, F is onto.

Second, we will show F preserves the metric. Let $\overline{\{A_{n}\}}$, $ \overline{\{B_{n}\}}$ $\in $  $\overline{\Omega}$. Then we have f(n) and g(n) as before.

\[
\begin{aligned}
 \mu (   A_{f(n)}  \triangle B_{g(n)}  ) = &  \mu^{\ast} (   A_{f(n)}  \triangle B_{g(n)}  )  \\
  \leq  &   \mu^{\ast} (   A_{f(n)}  \triangle  \overline{\lim} A_{f(n)} )  +   \mu^{\ast} (  \overline{\lim} A_{f(n)}  \triangle \overline{\lim} B_{g(n)}  )       +   \mu^{\ast} (  \overline{\lim} B_{g(n)}  \triangle B_{g(n)} )
 \end{aligned}
\]
By taking the limit on both sides,
\[
  \lim \mu (   A_{n}  \triangle B_{n}  )  =  \lim \mu (   A_{f(n)}  \triangle B_{g(n)}  ) \leq  \mu^{\ast} (  \overline{\lim} A_{f(n)}  \triangle \overline{\lim} B_{g(n)}  )   .
\]

In addition,
\[
\begin{aligned}
& \mu^{\ast} (  \overline{\lim} A_{f(n)}  \triangle \overline{\lim} B_{g(n)}  )   \\
&\leq     \mu^{\ast} (   A_{f(n)}  \triangle  \overline{\lim} A_{f(n)} )  +   \mu (   A_{f(n)}  \triangle B_{g(n)}  ) +  \mu^{\ast} (  \overline{\lim} B_{g(n)}  \triangle B_{g(n)} )  .
\end{aligned}
\]
By taking the limit on both sides,
\[ \mu^{\ast} (  \overline{\lim} A_{f(n)}  \triangle \overline{\lim} B_{g(n)}  )
\leq  \lim \mu (   A_{f(n)}  \triangle B_{g(n)}  ) .
\]
Therefore, $\overline{d}( \overline{\{A_{n}\}}, \overline{\{B_{n}\}} ) $ = $\lim \mu (   A_{n}  \triangle B_{n}  ) $ = $\mu^{\ast} (  \overline{\lim} A_{f(n)}  \triangle \overline{\lim} B_{g(n)}  )     $

=  $\mu^{\ast} (  F ( \overline{\{A_{n}\}} )   \triangle   F ( \overline{\{B_{n}\}} )  )     $ = $d(  F ( \overline{\{A_{n}\}} )  ,   F ( \overline{\{B_{n}\}} )     ) $.

Lastly, we will show that F is one to one. Let $  F ( \overline{\{A_{n}\}} )  ,   F ( \overline{\{B_{n}\}} )    \in \widetilde{\mathbf{S}}$ such that $  F ( \overline{\{A_{n}\}} )  =  F ( \overline{\{B_{n}\}} ) $ a.e..

Then $ \overline{\lim} A_{f(n)}  =  \overline{\lim} B_{g(n)}   $ a.e. implies $\mu^{\ast} (  \overline{\lim} A_{f(n)}  \triangle \overline{\lim} B_{g(n)}  ) = 0$. Then, as in the proof of F being onto,  $\lim \mu (   A_{n}  \triangle B_{n}  ) $ = $\mu^{\ast} (  \overline{\lim} A_{f(n)}  \triangle \overline{\lim} B_{g(n)}  )     $. Thus $\overline{\{A_{n}\}}$ = $\overline{\{B_{n}\}}$. Thus, F is one to one. Therefore, F is an isometry between $\overline{\Omega}$ and $ \widetilde{\mathbf{S}}_{\diagup_{\sim} }$.

\end{proof}
\end{theorem}

\begin{definition}\label{d:6}
Suppose X and Y are $\sigma$-algebras, and F: X $\rightarrow  $ Y is a one to one, onto well defined map. Then F is called a $\sigma$-algebra isomorphism if
\[
\begin{aligned}
F(\cdot  \bigcup \cdot ) = & F(\cdot ) \cup  F(\cdot )  , \  \  \ \  \ \  \  F(\bigcup^{\infty}_{i=1} \cdot ) & = &  \cup^{\infty}_{i=1}   F(\cdot )  ,  \\
F(\cdot  \bigcap \cdot ) = & F(\cdot ) \cap F(\cdot )  ,\ \ \  \ \ \ \ F(\cdot^{\textbf{C}} ) & = & F(\cdot )^{\textbf{C}} .
\end{aligned}
\]

\end{definition}

\begin{lemma}\label{le:1}

Let $E_{i}$ =  $\overline{\{B^{i}_{n}\}}$  $\in \overline{\Omega}$ for i $\geq$ 1 and by following Lemma 8 in [5], construct $Y_{L} $ =  $\cup_{i=1}^{N_{L}} B^{i}_{K_{L}}$ for each $L$ such that

\[
\mu^{\ast}(\cup_{i=1}^{\infty}   S_{i}  \triangle  \cup_{i=1}^{N_{L}} B^{i}_{K_{L}}     ) \leq  \mu^{\ast}(\cup_{i=N_{L}+1}^{\infty}   S_{i}) + \mu^{\ast}(\cup_{i=1}^{N_{L}}   S_{i}  \triangle  \cup_{i=1}^{N_{L}} B^{i}_{K_{L}}     ) < \frac{1}{L}.
\].

Then
$\overline{\{ Y_{L}   \}}  $ = $\bigcup_{i=1}^{\infty}   E_{i}  $.

\begin{proof}
Note that $E_{i}$ =  $\overline{\{B^{i}_{n}\}}$ = $\overline{\{B^{i}_{K_{L}}\}}$.

$( \bigcup_{i=1}^{n}   E_{i}  ) \bigcap \overline{\{ Y_{L}   \}} $  =  $ \overline{\{  \cup_{i=1}^{n} B^{i}_{K_{L}}   \cap  Y_{L}         \}}   $ = $ \overline{\{  \cup_{i=1}^{n} B^{i}_{K_{L}}     \}}   $ = $\bigcup_{i=1}^{n}   E_{i}   $ for any n.

Let $N_{L} > n$.
\[
\begin{aligned}
\mu( \cup_{i=1}^{N_{L}} B^{i}_{K_{L}}  \cap  (\cup_{i=1}^{n} B^{i}_{K_{L}})^{\textbf{C}}   ) = & \mu( \cup_{i=1}^{N_{L}} B^{i}_{K_{L}}  \triangle  \cup_{i=1}^{n} B^{i}_{K_{L}}     )  = \mu^{\ast}( \cup_{i=1}^{N_{L}} B^{i}_{K_{L}}  \triangle  \cup_{i=1}^{n} B^{i}_{K_{L}}     ) \\  \leq & \mu^{\ast}(\cup_{i=1}^{\infty}   S_{i}  \triangle  \cup_{i=1}^{N_{L}} B^{i}_{K_{L}}     ) + \mu^{\ast}(\cup_{i=1}^{\infty}   S_{i}  \triangle  \cup_{i=1}^{n} B^{i}_{K_{L}}     ).
\end{aligned}
\]

This implies that $\lim  \overline{\mu}( \overline{\{ Y_{L}   \}}  \bigcap ( \bigcup_{i=1}^{n}     E_{\alpha_{i}}) ^{\textbf{C}}  )   $ = 0. Therefore by the uniqueness of $\bigcup_{i=1}^{\infty}   E_{i}  $, $\overline{\{ Y_{L}   \}}  $ = $\bigcup_{i=1}^{\infty}   E_{i}  $.

\end{proof}

\end{lemma}

\begin{theorem}\label{t:4}
F is a $\sigma$-algebra isomorphism between $\overline{\Omega}$ and $ \widetilde{\mathbf{S}}_{\diagup_{\sim} }$.

\begin{proof}
We already showed that F is a one to one, onto map in Theorem $\ref{t:3}$.

Since, in general, $ \overline{\lim} A_{n}  \cup B_{n} $ = $ \overline{\lim} A_{n}  \cup \overline{\lim} B_{n} $, $F(\cdot  \bigcup \cdot )$ =  $F(\cdot ) \cup  F(\cdot )  $ follows immediately.

Let $\overline{\{B_{n}\}} \in \overline{\Omega}$.

Then,
\[
F ( \overline{\{B_{n}\}}^{\textbf{C}} )  =    F ( \overline{\{  ( B_{n})^{\textbf{C}}   \}} )   =                    \overline{\lim}( B_{f(n)} )^{\textbf{C}}  = (  \underline{\lim} B_{f(n)} )^{\textbf{C}}  = (  \overline{\lim} B_{f(n)} )^{\textbf{C}}  a.e..
\]
Note: by the construction of f(n), $\underline{\lim} B_{f(n)} = \overline{\lim} B_{f(n)}$ a.e. Thus, $F(\cdot^{\textbf{C}} )$ =  $F(\cdot )^{\textbf{C}} $ in $ \widetilde{\mathbf{S}}_{\diagup_{\sim} }$.

Similarly, $F(\cdot  \bigcap \cdot )$ = $F(\cdot^{\textbf{C}} \bigcup \cdot^{\textbf{C}} )^{\textbf{C}} $ =  $[ F(\cdot^{\textbf{C}} ) \cup  F(\cdot^{\textbf{C}} ) ]^{\textbf{C}} $= $F(\cdot ) \cap F(\cdot )  $.

Let $E_{\alpha_{i}}$ $\in \overline{\Omega}$ for i $\geq$ 1 and $E_{\alpha_{i}}$ = $\overline{\{B^{\alpha_{i}}_{n}\}}$.

Then for each i, there exists a $S_{i} = \overline{\lim} B^{\alpha_{i}}_{f(n)}  \in   \widetilde{\mathbf{S}}  $ such that $ \lim \mu^{\ast} (B^{\alpha_{i}}_{n} \triangle S_{i})  =0  $.

Now suppose we have $\{ Y_{L} \}$ in the same manner as Lemma $\ref{le:1}$. By design, $\{ Y_{L} \}$ converges to $\cup_{i=1}^{\infty}   S_{i} $. Then $\overline{ \{ Y_{L} \}  }$ = $\bigcup_{i=1}^{\infty}   E_{\alpha_{i}}  $ by Lemma $\ref{le:1}$. Now we have

\[
F( \bigcup_{i=1}^{\infty}   E_{\alpha_{i}}  )  = F(\bigcup_{i=1}^{\infty}   \overline{\{B^{\alpha_{i}}_{n}\}}  )  = F ( \overline{ \{ Y_{L} \}  } ) = \overline{\lim} Y_{f(L)}   .
\]
Since $\lim \mu^{\ast} ( Y_{L} \triangle \overline{\lim} Y_{f(L)}  )  =0  $ and $\lim \mu^{\ast} ( Y_{L} \triangle \cup^{\infty}_{i=1} S_{i}  )  =0  $, we have $\overline{\lim} Y_{f(L)}   $ = $\cup^{\infty}_{i=1} S_{i}  $ a.e.. In addition, $\cup^{\infty}_{i=1} S_{i}  $ = $\cup^{\infty}_{i=1} F(  \overline{\{B^{\alpha_{i}}_{n}\}}  )  $.

Thus,
\[
F( \bigcup_{i=1}^{\infty}   E_{\alpha_{i}}  )  = \cup^{\infty}_{i=1} F( E_{\alpha_{i}} )  .
\]
Therefore, the claim follows.

\end{proof}
\end{theorem}

\section{Conclusion}
Theorem $\ref{t:3}$ and  $\ref{t:4}$ show that the completion of $\Omega$ is isometric and $\sigma$-algebra isomorphic to $ \widetilde{\mathbf{S}}_{\diagup_{\sim} }$. Thus the completion of $\Omega$ is isometric and $\sigma$-algebra isomorphic to the Catheordory Extension under the equivalence relation $\sim$ by the conclusion in [5].

\section{Acknowledgement}
  I would like to thank my grandfather Waichi Tanaka for his inspiration and financial assistance and Andrew Aames for encouraging him to progress through the graduate program. With the kind support of both, I have progressed further than I ever thought possible. In addition, I would like to thank my friends Richard Han and Eli Depalma for their editing assistance and Vincent Davis, Aaron Hudson, Mark Tseselsky for representing me and for their professional advice.



\end{document}